\documentclass[12pt]{article}

\usepackage{amsmath}
\usepackage{amssymb}
\usepackage{amsfonts}
\usepackage{amsthm}
\usepackage{latexsym}
\usepackage{epsfig}
\usepackage{array}
\usepackage{boxedminipage}

\usepackage{dsfont}%vector of one's - $\mathds{1}$

\makeatletter
\renewcommand{\section}{\@startsection
{section} {1} {0mm} {-\baselineskip} {0.5\baselineskip}
{\large\bf}}
\renewcommand{\subsection}{\@startsection
{subsection} {2} {0mm} {-\baselineskip} {0.5\baselineskip}
{\normalsize\bf}} \makeatother

\newcommand{\R}{\mathbb{R}}

\newtheorem{theor}{Theorem}

\newtheorem{coro}{Corollary}

\newcommand{\rk}{\mathrm{rk} \hspace{1mm}}
\newcommand{\cork}{\mathrm{cork} \hspace{0.5mm}}
\newcommand{\ke}{\mathrm{ker} \hspace{0.5mm}}

\newcommand{\dsp}{\displaystyle}

\newcommand{\W}{W\hspace{-1mm}}

\newcommand{\1}{_{\bf 1}}
\newcommand{\2}{_{\bf 2}}
\newcommand{\ml}{\mathbb{L}}

\newcommand{\mg}{\mathbb{G}}
\newcommand{\tra}{{\sf T}}

\begin{document}

%\maketitle

\begin{center}

{\Large\bf Some qualitative problems in network dynamics\footnote{Supported by Research Project MTM2010-15102
of Mi\-nis\-terio de Ciencia e Innovaci\'{o}n, Spain.}}

\vspace{4mm}

{\sc 
Ricardo Riaza}\\ %\footnote{Corresponding author}}\\
\ \vspace{-3mm} \\
Depto.\ de %Departamento de
Matem\'{a}tica Aplicada 
a las TIC \\
%a las Tecnolog\'{\i}as de la Informaci\'{o}n y las Comunicaciones\\
%Escuela T\'{e}cnica Superior de Ingenieros de 
ETS Ingenieros de  
%ETSI 
Telecomunicaci\'{o}n  \\ %, % \\
Universidad Polit\'{e}cnica de Madrid \\  %\\ Ciudad Universitaria s/n 
28040 Madrid, Spain \\
{\sl ricardo.riaza@upm.es} \\

%\vspace{0.5cm}

%[DRAFT]

\end{center}

\vspace{2mm}

\begin{abstract}
This paper addresses analytical aspects of deterministic, continuous-time
dynamical systems defined on networks. 
The goal is to model and analyze certain phenomena which must be framed
beyond the context of networked dynamical systems, understood 
as a set of interdependent dynamical systems defined on the nodes of a (possibly 
evolving)
%weighted and/or directed) 
graph. %network.
%which is a common working scenario in this area. 
In order to advocate for a more flexible approach to the study of network
dynamics, we tackle some qualitative problems which do not fall in 
this working scenario. First, we address a stability problem
on a network of heterogeneous agents, some of which are of dynamic nature 
while others just impose restrictions on the system behavior.
%; this naturally leads to a differential-algebraic formalism.
Our second context assumes that the network is clustered, and we address
a two-level stability problem involving the dynamics of both individual
agents and groups. The aforementioned problems exhibit lines of non-isolated equilibria,
and the analysis implicitly assumes a positiveness
condition on the edge weights; the removal of this restriction complicates matters, and
our third problem concerns a graph-theoretic characterization of 
the equilibrium set in the dynamics of certain networks with positive and negative
weights, the results applying in particular to signed graphs.
Our approach combines graph %(and hypergraph) 
theory
and dynamical systems theory, 
but also uses specific tools coming from linear algebra,
including e.g.\ 
Ger\v{s}gorin discs or Maxwell-type determinantal expansions.
The results are of application in social, economic, and flow networks, among others, and
are also aimed at motivating further research.
% and advocate
%a more flexible approach to the analysis of network dynamics.

%OLD (up to Jan 2014): oriented to ``constrained dynamics on networks''...:
%In a continuous-time, deterministic, lumped setting,
%the main attention is directed to the presence of constraints
%leading to a differential-algebraic formalism.
%Constraints arise naturally e.g.\ under an equilibrium
%assumption on the node dynamics in coevolution networks,
%or in problems in which --heterogeneous-- nodes --of different nature-- coexist
%(ie some nodes accommodate a state variable while others do not). 
%In this context we characterize the behavior of a semistate version
%of the dynamical system defined by the Laplacian matrix.
%Flows...
%Higher index DAEs arise for instance in the presence of
%constraints involving a given set of network objects (e.g.\ nodes)
%in problems in which the model involves (nmg) others (e.g.\ links),
%or in multilevel networks in which certain variables are only visible
%at a given level.

%ALGO MAS, EJEMPLOS, ETC.

\end{abstract}

\vspace{2mm}

\noindent {\bf Keywords:} graph; digraph; signed graph; hypergraph;
network; dynamical system; stability; equilibrium; 
differential-algebraic equation; flow; social network; economic network.

\

%\vspace{5mm}

\noindent {\bf AMS subject classification:} %2010
05C21, %Flows in graphs
05C22, %Signed and weighted graphs
05C50, %Graphs and linear algebra
%05C70, %Factorization, matching, partitioning, covering and packing
05C82, %Small world graphs, complex networks
%15A22, %Matrix pencils
%34A09, %Implicit equations, differential-algebraic equations [See also 65L80]
%34C23, %Bifurcation (ODEs)
%34C45, %   	Invariant manifolds
34D20, %ODEs: stab. theory, Lyap. stab.
%34D35, %   	Stability of manifolds of solutions
%37C75, %Smooth dyn. sys., Stab. theory
%37G10, %Bifurcation of singular points;
90B10, %Network models, deterministic
91D30, %Social networks
%94C05. %Analytic circuit theory
%94C12 Fault detection; testing
94C15. %Applications of graph theory

%\newpage

%\setcounter{page}{1}

\section{Introduction}
\label{sec-intr}

%Somewhere multidisciplinary field... first paragraph

The interactions between dynamical systems and networks have attracted
the interest of many researches in the last decade. Needless to say,
the study of 
dynamics taking place on graph-theoretic structures 
has a long tradition in different application fields, including
nonlinear circuits, water and gas networks, 
neural nets, 
%dynamic flows (transit times en FF type
%``dynamic flows in networks'' or ``flows in dynamic networks'');
flows, coupled dynamical systems and, in particular, 
coupled oscillators and synchronization, 
%(\cite{kaob, strogatzreview, zelazo}, MORE);
etc., to name but a few. 
%Remove: biological sciences; traffic;
%transportation networks (Peruani en Dyn on and of...): entities at nodes, part individuals, migration, diffusion rate for the flowrate, nmg (confuso)
In the framework of what is now known
as {\em network science,} dynamics is one of the key research directions. 
A great deal of research within this
interdisciplinary field is directed
to the dynamics of state variables defined on the nodes of a network;
the topics involved include networks of coupled dynamical systems 
\cite{kaob, newmannetworks, strogatzreview, zelazo},
consensus and agreement protocols
\cite{%olfati2003, 
olfati2007, rahmani2009, tanner2004, yangliu},
controllability and related aspects \cite{liu2011, liu2013, rahmani2009}, 
%Pello, Arkak? Bergson arxiv 2011 cited Pello;
epidemics \cite{easley, newmannetworks},
diffusion of information, innovations, etc.\ \cite{barrat, mmsna},
%iterated games \cite{newmanreview} y refs 1, 135, 230, 415 alli;
%solution via eigenvectors of the Laplacian, motiva esto)
dynamics of social networks \cite{mmsna, holmenewman, lozano}
%so-called longitudinal social network data \cite{mmsna};
or evolutionary dynamics \cite{albertbarabasireview,
lieberman}, among others. %OK
%phase transitions on networks \cite{newmanreview}, %Ising
%VER EL DE HOLME Y NEWMAN ARRIBA
%Synchronization in small-
%world systems: Barahona, M. \& Pecora, L. M. 2002 Synchronization in small-
%world systems. Phys. Rev. Lett. 89, 054101-4. (doi:10.
%1103/PhysRevLett.89.054101)
%convergence of opinion in social systems \cite{holmenewman},
%Controllability of complex networks and related aspects have also attracted
%attention from different perspectives: 
%cf.\ \cite{liu2011, liu2013, rahmani} %Pello, Arkak?
%and references therein.

Much attention has also been paid
to the study of 
network growth  \cite{albertbarabasireview, barabasi99,
price65, price76}, %Price: network of citations
percolation \cite{newmanreview}
and, more generally, to the analysis of the
dynamic changes of the network itself
%Other references in this direction are 
(cf.\ \cite{albertbarabasireview, newmannetworks, siljak}). 
%percolation (=filtrado), modelled via removing of nodes, e.g.
%failure of routers or vacunado 
These dynamic changes of the network structure have been studied
in the last decades both in a deterministic and in a probabilistic
context.
The original ideas regarding these
{\em evolving networks}
\cite{albertbarabasireview} can be traced back to 
the seminal work of Erd\H{o}s and R\'enyi on the evolution of
random graphs \cite{erdos1959, erdos1960}. Find more references
in \cite{albertbarabasireview, newmanreview, newmannetworks}.
In a deterministic, continuous-time context,
the changes in the structure of a network modelled
as a weighted simple graph can been naturally described
as a dynamical system in which the state variables are 
the graph weights or, equivalently, 
the adjacency matrix, in the understanding
that a vanishing weight $w_{ij}$ models the absence of
an edge between nodes $i$ and $j$. This dynamical system
can be defined as a one-parameter group of transformations
of a given state space,
%*** check the name, automosphisms, homeomorphisms ***
%as it is customary in topological dynamics, 
or as a system of differential equations. Find details
in \cite{siljak}. 

So-called {\em coevolutionary} processes involve simultaneously 
dynamics taking place on a network and changes in the network structure;
cf.\ \cite{fronczak, gross, holmenewman, koenigbattiston,lozano} 
and references therein.
Indeed, it is natural to consider networks in which dynamical changes
in the nodes coexist with the evolution of the network
itself. With different names,
these coevolutionary processes have been addressed
 in several application fields, 
including neural network training
\cite{baldi, haykin, nnspp}, %Tang y Siljak 
social network analysis \cite{holmenewman, lozano},
economic models \cite{koenigbattiston}, etc.
%Tb K\"onig, economic models, 
%distinguish between the dynamics taking place on the
%agents (nodes) and the evolution of the network. Also different
%time scales.
See also \cite{fronczak, gross} and the bibliography therein.
%IMP: Gross \cite{gross}
%explicitly terms this ``coevolution'' of dynamics on the NODES
%and dynamics of the network itself.

%Mas COevolution (ojo, todo lo saco de gross):
%Gil, S. \& Zanette, D. H. 2006 Coevolution of agents and
%networks: opinion spreading and community disconnection. Phys. Lett. A 356, 89–95. (doi:10.1016/j.physleta.
%2006.03.037)
%Zimmermann, M. G., Egu\'iluz, V. M. \& San Miguel, M. 2004
%Coevolution of dynamical states and interactions in
%dynamic networks. Phys. Rev. E 96, 065102-4. (doi:10.
%1103/PhysRevE.69.065102)

%THIS I REVOMVE, DEC 2014:

%Less attention has been paid to {\em constrained} dynamics... leading to a 
%DAE formalism...

%Exceptions: 

%1) Multiagent described by DAEs \cite{yangliu}; 

%2) Lennart. 

%3) Neural learning under an equilibrium assumption (Baldi, spp) 

%Llevo a final seccion 2:
%4) Potential-driven flow dynamics in tree networks %(ok)  
%\cite{mayes2011} {mayes2011} 

%5) Integrated process networks %(ok) 
% \cite{contou2004} {contou2004} 

%6) Pressurized water networks %(ok)
%\cite{steinbach2005} {steinbach2005} 

%7) Dynamic traffic networks %ok
%\cite{han2012} {han2012} 

%\section{Dynamics on networks}
%\label{sec-dyn}

In spite of the large amount of literature addressing dynamical processes on
networks,
it seems that in a continuous-time, deterministic context, 
general approaches to the formulation and
analysis of network dynamics are somehow too restrictive, being
focused on one of just three paradigms: %the following three settings:
a) systems in which the dynamics takes places only on the network nodes; b)
dynamics of the network itself; and c) coevolution processes
combining both. The references cited above support this point of view. Certainly,
in specific fields, dynamical systems on networks arise beyond
these categories; flow networks or nonlinear electrical circuits, which systematically 
exhibit state variables not supported on the network nodes, are examples of this. But a general
approach, studying abstract dynamics on networks beyond the three settings indicated above, 
seems to be lacking in the literature.

This paper is intended to be a step in this direction, addressing certain dynamical
problems which cannot be framed in the working scenarios mentioned above.
First, models for network dynamics formulated in terms of a set of dynamical systems defined on the nodes of a graph somehow
assume a homogeneous nature on all agents. By contrast, 
we address in Section \ref{sec-constrained} a stability problem involving
a potential-driven flow
on a network of heterogeneous agents, some of which are of dynamic nature 
while others just impose restrictions on the system behavior; 
this naturally leads to a differential-algebraic formalism.
In Section \ref{sec-groups} we will examine clustered networks, addressing
a two-level stability problem  which involves the dynamics of both individual
agents and groups. 
Lines of non-isolated equilibria will naturally
arise in the two problems just mentioned; the one-dimensional nature of the
equilibrium set is implicitly related to
a positiveness condition on the edge weights. 
In (possibly nonlinear) flow problems on networks
which may include negative weights, the structure of 
the equilibrium set becomes more involved: 
in Section \ref{sec-eq}
we will examine such structure without the aforementioned positiveness 
assumption, applying in particular
the results to flow networks on signed graphs. Finally, Section
\ref{sec-con} compiles some concluding remarks.

%\subsection{Dynamics on nodes}
%\label{subsec-nodes}

%Newman; consensus; NDS de Zelazo;
%coupled oscillators; Regular networks of coupled dynamical systems 
%(Strogatz 2001);
%CSNs (Coupled Systems on Networks): \cite{kaob}

%[Zelazo: ``networked dynamical systems''. Linear. Idea: sistemas
%dinamicos en los agentes, acoplados por entrada o salida o estado o combinaciones. Ejemplo que da de coupled by state: consensus.]

%\subsection{Dynamics of networks and coevolution}
%\label{subsec-of}

%%%%%%%%%% LLEVO A SECC 1
%Dynamic changes of the network structure have been studied
%in the last decades both in a deterministic and in a probabilistic context.

\section{Heterogeneous agents and constrained dynamics}  %on digraph-based networks}
\label{sec-constrained}

Along the terms presented in the Introduction, the ``dynamics on nodes'' scenario 
for modelling continuous-time, deterministic 
network dynamics assumes that the dynamical process can be described by 
an explicit ODE such as
\begin{equation} \label{gensys}
x'=f(x),
\end{equation}
where $x$ is the state vector, which
is defined by a set of variables supported on the network nodes. 
The topology of the network is assumed to be comprised on the vector field
$f$ and, for the sake of simplicity, we are ignoring control variables,
explicit time dependences associated with system inputs, etc.
A very common perspective assumes
that the network edges model a coupling between a set of (otherwise
independent) dynamical systems defined on the network nodes
(see e.g.\ \cite{kaob, newmannetworks, strogatzreview, zelazo}), in a
way such that only the states of adjacent systems are visible by each node.
This means that the differential equations above read
as
\begin{equation} \label{local}
x_i' = f_i(x_i, x_{j_1}, \ldots, x_{j_{d(i)}}), 
\end{equation}
for the $i$-th node ($i=1, \ldots n$); here $x_i$ is the set of state variables defined on node
$i$, and $j_1, \ldots, j_{d(i)}$ are the nodes adjacent to $i$.
Note that for a fully connected (complete) network, (\ref{local}) makes no
difference to (\ref{gensys}) (this is the case for instance in the
Kuramoto model; see e.g.\ \cite{strogatzreview}).

%\subsection{Constrained dynamics on networks}

Needless to say, the setting described above has a very broad application scope,
as the literature discussed in Section \ref{sec-intr} shows. However, in other 
cases this may be
too restrictive a framework to model certain dynamical processes on networks.
In particular, systems (\ref{gensys}) and (\ref{local}) do not capture the
eventual existence of constraints restricting the 
%admissible set of values for the state variables.
whole dynamics of the network. These constrained dynamical systems
are better framed in the theory of 
differential-algebraic equations (DAEs; cf.\ \cite{bre1, kmbook, LMTbook,
wsbook}).

%*** Motivate constrained etc. CONECTANDO CON LO DE ABAJO
%above are homogeneous, although the fi's may be different
%Constrained formulations  may arise e.g.
%considering a network with dynamical processes taking place 
%on the nodes, and in which these node dynamics are 
%themselves modelled as a constrained system. 
%See e.g.\ the consensus protocols for DAE agents discussed in
%\cite{yangliu}. This is an interesting topic but
%not in the scope of the present paper.
%Broadly speaking, in this section and in the next one 
%we will focus instead 
%on problems in which explicit constraints restrict the 
%whole dynamics of the network.
%For ease of discussion
%To fix ideas, we will consider networks 
%modeled by possibly weighted (di)graphs 
%in which the state (or, better, semistate) variables are
%defined on the nodes (denoted by $x$), 
%on the edges (denoted by $u$), or are the weights $w$ themselves;
%that is, we disregard the variable $z$ and the equation
%for $z'$ in  (\ref{expl}).
%(CAREFUL USE Z FOR ANOTHER THING, ``OTRAS VARIABLES'')

In this context, we focus in this section on networks in which
%First, it may well happen that in a network 
two different sets of nodes (agents) coexist; the behavior of some of them 
are defined by explicit dynamical systems, whereas the others
do not involve dynamics explicitly
but are    %in the state variables and the other 
governed, instead, by algebraic (non-differential) restrictions.
%If there are no edge variables and the weights are constant, 
Letting
$y$ stand for the dynamic variables, which model the states of agents with
explicit dynamics, and $z$ for the algebraic variables, which correspond
to agents without explicit dynamics,
these systems can be written as a so-called semiexplicit DAE
\begin{subequations}\label{dae2}
\begin{eqnarray} 
y' &= & f(y, z) \\
0 & = & g(y, z), 
\end{eqnarray}
\end{subequations}
where $f$ and $g$ capture the topology of the network; an example can be found in (\ref{consDAE}) below.
%This will the case for partial storage subsection \ref{subsec-het}
%Even if all agents have explicit dynamics, the nature of the system
%may impose a set of external constraints. Disregarding weights, edge variables, etc. 
%a constrained autonomous system may take the form
%\begin{subequations}
%\begin{eqnarray} 
%x' &= & f(x) \\
%0 & = & g(x). 
%\end{eqnarray}
%\end{subequations}
%Provided that $x \in \R^n$, $f:\R^n \to \R^n$, this is in general an overdetermined system... (refs). 
%Instead, in many cases it is realistic to assume that some other variables
%(e.g.\ some edge variables $u$) are not explicitly 
%fixed by the states $x$ of the nodes;
%%, and that they fluctuate with the dynamics; 
%in this setting the model reads, instead,

%where $x$... $u$... $f$... $g$...
%This is a Hessenberg DAE (REFS)...
%In the model (\ref{dae1}), equation (\ref{dae1d}) may
%well model a set of external constraints, imposed by the nature
%of the problem. In particular, if the constraints involve
%only one type of variables (e.g.\ state variables defined on the nodes)
%in systems with more one than type (e.g.\ when edge variables
%are also concerned) we may well get a higher index DAE;
%cf.\ hess 2 in Sect. 3.3... (and refer to groups Sect 4)
%An example of this will be discussed in subsection \ref{subsec-hess2}

%when different types of actors (nodes) coexist. We focus on a problem in which
%some of the actors have a dynamic description whereas others do not involve dynamics explicitly.

\subsection{Potential-driven flows and the graph Laplacian}
\label{subsec-flows}

To fix ideas, consider a simple connected network in which a given  
%In order to do so we begin with an alternative interpretation of the 
%dynamical system (\ref{consensus}) associated with the graph Laplacian. We assume
%the graph to be simple, and will look at
%(\ref{consensus}) as a process in which the actors of a network collect and
%redistribute a given 
resource or commodity flows among a set of $n$ agents located at the network nodes; here
``simple'' means that the network has neither parallel edges nor self-loops.
Let us first assume that all
agents collect this commodity, and denote the amount of 
the collected resource at node $i$ by $x_i$.  
The flowrate (or flow) of this resource at a given edge $j$, connecting 
nodes $i$ and $k$, will be denoted by $u_j$. Every edge is endowed with a reference 
direction; if the $j$-th edge is directed 
say from $i$ to $k$, $u_j>0$ (resp.\ $u_j <0$) means that the resource flows from 
$i$ to $k$ (resp.\ from $k$ to $i$). This drives the problem to the
context of directed graphs, 
although it is worth emphasizing that all the results are
independent of the choice of directions.
Denote $x=(x_1, \ldots, x_n)$, $u=(u_1, \ldots, u_m)$,
where $m$ is the number of edges.

With these reference directions, the entries of the incidence matrix
$A = (a_{ij}) \in \mathbb{R}^{n \times m}$ of the resulting digraph read as
\begin{eqnarray} \label{incidence}
a_{ij} = \left\{
\begin{array}{rl}
1 & \text{ if edge } j \text{ leaves  node } i \\
-1 & \text{ if edge } j \text{ enters node } i \\
0 & \text{ if edge } j \text{ is not incident with node } i.
\end{array}
\right.
\end{eqnarray}
The continuity  equations at 
the nodes can be then written as
%\begin{subequations} \label{unfold-lapl}
\begin{eqnarray}
x' & = & -Au, \label{unfold-lapla} %\\
%u & = & A^{\tra} x. \label{unfold-laplb}
\end{eqnarray}
%\end{subequations}
since the increase $x_i'$ of the stored resource at node $i$  equals 
the net flow
entering the node, that is, $-A_i u$ with the sign convention above and $A_i$ standing
for the $i$-th row of $A$.

System (\ref{unfold-lapla}) is underdetermined. The way in which the flowrates $u$
interact with the node state variables $x$ may be defined from different perspectives; 
for instance,
in a game-theoretic setting the agents' strategies would define the flows, or 
in the framework of control theory $u$ might be designed as to achieve a given goal.
Here we will simply assume that $u$ is designed in order to get a fair (equal) distribution of
the commodity among all agents. This can be achieved by setting the following potential-driven
definition of the flowrates:
\begin{eqnarray}
%x' & = & -Au \label{unfold-lapla} \\
u & = & A^{\tra} x, \label{unfold-laplb}
\end{eqnarray}
which is just a redistribution law 
in which the flow from $i$ to $k$ equals $x_i - x_k$. It is a trivial matter to recast
(\ref{unfold-lapla})-(\ref{unfold-laplb}) as
\begin{equation}
x' = -\mathbb{L}x \label{consensus}
\end{equation}
where $\mathbb{L}=AA^{\tra}$  is the graph Laplacian matrix. Note that (\ref{consensus})
has the form depicted in (\ref{local}), as a consequence of
the identity $\mathbb{L} = -C+D$, where $C$ is the adjacency 
matrix of the graph (without directions) 
and $D$ is a diagonal matrix whose $i$-th
entry is the degree of node $i$.
%Now, by using the identity $\mathbb{L}=AA^{\tra}$ and gluing together the $x_i$'s and $u_k$'s
%into vectors $x$, $u$,
%we may unfold (\ref{consensus}) as
%\begin{subequations} \label{unfold-lapl}
%\begin{eqnarray}
%x' & = & -Au \label{unfold-lapla} \\
%u & = & A^{\tra} x. \label{unfold-laplb}
%\end{eqnarray}
%\end{subequations}
%Equation (\ref{unfold-lapla}) is nothing but the continuity equations at 
%the nodes, since the increase $x_i'$ of the stored resource at node $i$  equals the net flow
%entering the node, that is, $-A_i u$ with the sign convention above. 
%In turn, (\ref{unfold-laplb}) can be seen as a
%redistribution law, in which the flow from $i$ to $j$ equals $x_i - x_j$. Mathematically, 
%system
%(\ref{unfold-lapl}) simply presents the graph Laplacian (or consensus) dynamics
%as a system of continuity equations in which the flowrates 
%are potential-driven. As indicated above, the dynamics converge to ***
The dynamics of (\ref{consensus}) has been analyzed in the context of so-called
{\em consensus} protocols \cite{%olfati2003, 
olfati2007, rahmani2009, tanner2004, yangliu},
and the state variables $x$ may be checked to converge to a common value which is
the arithmetic mean of the initial values
%the trajectory of (\ref{consensus}) emanating
%from a given $x(0)=(
$x_1(0), \ldots, x_n(0).$ %)$ converges
%to a common value $$x_1= \ldots = x_n = 
%\mathrm{Average}[x_1(0), \ldots, x_n(0)].$$
%Different extensions and related results can be found in
%the above-cited references \cite{olfati2003, olfati2007, rahmani2009, tanner2004, yangliu}.
This means that the resource is indeed redistributed among all the agents
in a way such that all of them asymptotically store the same quantity.

\subsection{Heterogeneous agents}

We now drive our attention to a different problem; given again a simple 
connected network, assume that only some agents,
say type-1 agents (which w.l.o.g.\ are assumed to be those
from 1 to $r$, with $1 \leq r < n$) 
accumulate the aforementioned resource, again in an amount $x_i$ for the $i$-th agent. 
The continuity equations for type-1 agents read
as above,
\begin{eqnarray}
x_i'  =  -A_iu, \ i=1, \ldots, r. \label{equat1}
\end{eqnarray}
The other agents, say type-2 ones (those numbered 
from $r+1$ to $n$) are different and do not collect the commodity; instead, they are 
however intermediaries among type-1 agents.
Therefore, the continuity equations at type-2 agents are
\begin{eqnarray}
0  =  A_iu, \ i=r+1, \ldots, n. \label{equat2}
\end{eqnarray}
%Mention mixture of the cases for $\alpha$ in Definition ***
Finally, we assume that the flowrates are still given by (\ref{unfold-laplb});
note that by means of the variables $x_i$ type-2 agents now simply set up a reference 
value for the flowrates $u$. 
Denote $y = (x_1, \ldots, x_r)$, $z=(x_{r+1}, \ldots, x_n)$
(and, for later use, $y_i=x_i, \ z_i = x_{i+r}$) 
and split $A$ as
$$A = \left( \begin{array}{c} A\1 \\ A\2 \end{array} \right),$$
where $A\1$ (resp.\ $A\2$) is the submatrix of $A$ defined by the
first $r$ rows (resp.\ last $n-r$ rows).
The boldface subscripts in $A\1$ and $A\2$ 
are meant to distinguish these submatrices from the first and second row
of $A$, to be written as $A_1$, $A_2$.
%equations 
%(\ref{equat1}), (\ref{equat2}) and (\ref{unfold-laplb})
With this notation, the network dynamics is modelled by the DAE
\begin{subequations}
\begin{eqnarray}
y' & = & -A\1A^{\tra}x \\
0 & = & A\2 A^{\tra} x.
\end{eqnarray}
\end{subequations}
Splitting $x=(y,z)$ in the right-hand side, and
denoting $\ml_{ij}=A_{\bf i}A_{\bf j}^{\tra}$ for ease of notation, this DAE can
be rewritten as
%\begin{subequations}
%\begin{eqnarray}
%y' & = & -A\1A\1^{\tra}y - A\1A\2^{\tra}z \\
%0 & = & A\2 A\1^{\tra}y + A\2A\2^{\tra}z
%\end{eqnarray}
%\end{subequations}
%(compare with (\ref{consensus})).
\begin{subequations}\label{consDAE}
\begin{eqnarray}
y' & = & -\ml_{11}y - \ml_{12}z \\
0 & = & \ml_{21}y + \ml_{22}z. \label{consDAEb}
\end{eqnarray}
\end{subequations}
Note that this constrained system has the differential-algebraic
form displayed in (\ref{dae2}).
%is a semistate version of the Laplacian dynamics. 

Theorem \ref{th-constrained} describes the dynamics
of the DAE (\ref{consDAE}); in particular, it shows that
the scheme defined above redistributes the commodity
in a way such that all type-1 agents asymptotically accumulate the same amount 
of the resource. Note also that there is a somewhat unusual 
phenomenon, namely the existence of a line of non-isolated 
equilibria. We use from DAE theory the notion of a 
{\em consistent initial value}, which is a value of $(y(0), z(0))$ which
satisfies the constraints (\ref{consDAEb}) (cf.\ \cite{bre1, kmbook, LMTbook,
wsbook}).

\begin{theor} \label{th-constrained} 
In the setting described above,
the following assertions hold for the DAE (\ref{consDAE}).

\vspace{3mm}

\noindent {\rm (a)} The constraint (\ref{consDAEb}) specifies an $r$-dimensional
solution space. Its (transversal) intersection with the hyperplanes
$\sum_{i=1}^r y_i=k \in \R$ defines a foliation of the dynamics
by a family of $(r-1)$-dimensional
invariant spaces.

\vspace{3mm}

\noindent {\rm (b)} All equilibrium points are located
in the line
$y_1 = y_2 = \ldots = y_r = z_1 = \ldots = z_{n-r}$, 
which is transversal to the aforementioned
invariant spaces if $r \geq 2$.

\vspace{3mm}

\noindent {\rm (c)}  The trajectory $(y(t),z(t))$ emanating from a given consistent
initial value $(y(0),z(0))$
converges exponentially to the equilibrium
$$y_i = z_j = \displaystyle\frac{y_1(0)+ \ldots + y_r(0)}{r}, 
\ i = 1, \ldots, r, \ j=1, \ldots n-r.$$
\end{theor}

\

\noindent {\bf Proof.} 

\vspace{3mm}

\noindent (a) A well-known property in digraph theory states that any proper subset
of the rows of the incidence matrix $A$ of a connected digraph has maximal rank.
This implies that
the matrix $A\2$, defined by the last $n-r$ rows of $A$, has maximal row rank
(note that $r \geq 1$). In turn, this means
that $\ml_{22}=A\2 A\2^{\tra}$ is a non-singular matrix; indeed, provided that $A\2 A\2^{\tra}v=0$
we get $v^{\tra}\hspace{-1mm}A\2 A\2^{\tra}v=0$, and then $A\2^{\tra}v=0$ yields $v=0$.
Since $A\2 A\2^{\tra}$ has order $n-r$, it follows that (\ref{consDAEb})
specifies an $r$-dimensional linear space. 

This space is filled by solutions of the DAE, which are defined from those of
the explicit ODE
\begin{eqnarray}\label{reduced}
y'  =  -(\ml_{11} - \ml_{12}\ml_{22}^{-1}\ml_{21})y 
\end{eqnarray}
by means of the additional relation
\begin{eqnarray} \label{algeb}
z  =  -\ml_{22}^{-1}\ml_{21}y.
\end{eqnarray}

That (\ref{consDAEb}) intersects transversally the hyperplanes
$\sum_{1}^r y_i=k \in \R$ can be seen again as a consequence of the fact that $\ml_{22}=A\2 A\2^{\tra}$
is non-singular, since the coefficient matrix of the set of linear equations defining
the intersection reads as
$$
\left(\begin{array}{cc}
\mathds{1}_r^{\tra} & 0 \\
\ml_{21} & \ml_{22} 
\end{array}\right),
$$
where $\mathds{1}_r^{\tra} \in \R^{1 \times r}$ has all entries equal to one. 
This matrix is easily seen to have maximal row rank, meaning that
the intersection is indeed transversal and hence defines an $(r-1)$-dimensional space, 
for any fixed $k$.

Denoting by $\mathds{1}$ the vector of 1's in $\R^{n \times 1}$,
the invariance of these spaces follows from the identity $\mathds{1}^{\tra} A=0$ (expressing
the fact that the sum of all rows of the incidence matrix $A$ vanishes), which readily 
yields $\mathds{1}^{\tra} \ml=0$ since $\ml = A A^{\tra}$. Recasting (\ref{consDAEb})
as $-\ml x =0$, from (\ref{consDAE}) one gets $\mathds{1}_r^{\tra} y' = 0$,
meaning that $y_1 + \ldots + y_r$ indeed remains constant along trajectories.

\

\noindent (b) The equations $y_1 = y_2 = \ldots = y_r = z_1 = \ldots = z_{n-r}$ define a line
of equilibria because of the identity 
$A^{\tra} \mathds{1} =0$, which means that $\mathds{1}$ spans $\ker \ml
%yields $\ml\mathds{1} =0$ 
=\ker A A^{\tra}$. Therefore any vector belonging to span $<\hspace{-1mm}\mathds{1}\hspace{-1mm}>$ annihilates 
the right-hand side of (\ref{consDAE}) and therefore defines an equilibrium point of the DAE. 
There are no other equilibria because $\rk \ml = n-1$.

Note that if $r=1$, each one of the invariant spaces referred to in (a) amounts to a point, which
actually belongs to the equilibrium line $y_1=z_1=z_2=\ldots=z_{n-1}$. In cases with $r \geq 2$,
the transversality of the invariant spaces mentioned above and the equilibrium line
follows easily from the fact that, in $\R^n$, the hyperplane $\sum_{1}^r y_i=k$ intersects
transversally the line $y_1 = \ldots = y_r = z_1 = \ldots = z_{n-r}$. For later use
it is worth emphasizing
that the transversality of the intersection of the invariant spaces and the equilibrium line
implies that each equilibrium is unique within each one of such invariant spaces.

\

\noindent (c) The dynamics can be examined in $\R^r$ in terms of the $y$-variables via (\ref{reduced}),
having in mind that the $z$-components are given by (\ref{algeb}). In particular, the invariant spaces and the
line of equilibria discussed in (a) and (b) are projected respectively
onto the hyperplanes of $\R^r$ defined by
\begin{eqnarray} \label{invred}
\sum_{i=1}^r y_i=k
\end{eqnarray}
and the line span $<\hspace{-1mm}\mathds{1}_r\hspace{-1mm}>$, that is
\begin{eqnarray} \label{eqred}
y_1 = y_2 = \ldots = y_r.
\end{eqnarray}
Notice that both spaces are orthogonal to each other.

The rest of the proof relies on the fact that $-(\ml_{11} - \ml_{12}\ml_{22}^{-1}\ml_{21})$ is a
symmetric, negative semidefinite, corank-one matrix.
Its symmetry follows from that of $\ml_{11}$ and $\ml_{22}$, together with 
the identity $\ml_{21}=\ml_{12}^{\tra}$ (recall that $\ml_{ij}=A_{\bf i} A_{\bf j}^{\tra}$).
To check that it is negative semidefinite, write
$$ -v^{\tra}( \ml_{11} - \ml_{12}\ml_{22}^{-1}\ml_{21})v =
-\left(v^{\tra} \ \ -(\ml_{22}^{-1}\ml_{21}v)^{\tra}\right) 
\left(\begin{array}{cc}
\ml_{11} & \ml_{12} \\
\ml_{21} & \ml_{22} 
\end{array}\right) 
\left(\begin{array}{c}
v \\
-\ml_{22}^{-1}\ml_{21} v
\end{array}\right) \leq 0 
$$
for any $v \in \R^r$,
since $-\ml = -AA^{\tra}$ is itself negative semidefinite. Similarly, the fact that
$-(\ml_{11} - \ml_{12}\ml_{22}^{-1}\ml_{21})$ has corank one follows easily from
the remark that
$$ v \in \ke (\ml_{11} - \ml_{12}\ml_{22}^{-1}\ml_{21}) \Leftrightarrow 
\left(\begin{array}{c}
v \\
-\ml_{22}^{-1}\ml_{21} v
\end{array}\right) \in \ke \left(\begin{array}{cc}
\ml_{11} & \ml_{12} \\
\ml_{21} & \ml_{22} 
\end{array}\right) 
$$
together with the identity $\cork \ml = 1$.

Altogether, this implies that $-(\ml_{11} - \ml_{12}\ml_{22}^{-1}\ml_{21})$ has a
simple zero eigenvalue, its associated eigenspace being the line of equilibria
of (\ref{reduced}), and $r-1$ real and negative eigenvalues; the corresponding eigenvectors
are orthogonal to the equilibrium line because $-(\ml_{11} - \ml_{12}\ml_{22}^{-1}\ml_{21})$
is symmetric. Since the hyperplanes defined by (\ref{invred}) are themselves
orthogonal to the equilibrium line, the evolution in such invariant hyperplanes
is characterized by these $r-1$ negative eigenvalues. 

This means
that trajectories evolve exponentially towards the unique equilibrium in each invariant
hyperplane; it is a trivial matter to check that this equilibrium must be defined by 
$$y_i = \displaystyle\frac{y_1(0)+ \ldots + y_r(0)}{r}, \ i =1, \dots, r,$$
which is the unique solution to (\ref{invred})-(\ref{eqred}) with $k=y_1(0)+ \ldots + y_r(0)$.
Finally, the identities 
$$z_j = \displaystyle\frac{y_1(0)+ \ldots + y_r(0)}{r}, \ j =1, \dots, n-r,$$
follow from the fact that the equilibrium line spanned by $\mathds{1}_r$ in $\R^r$ is the
$y$-projection of the line spanned by $\mathds{1}$ in $\R^n$. \hfill $\Box$

\subsection{Constrained dynamics on general networks}

The systematic analysis of constrained dynamics on general networks defines a topic
of potential interest in many application fields. Just to name a few,
the semiexplicit DAE (\ref{dae2}) may,
in neural learning processes, 
model an equilibrium assumption for the dynamics of the neural state variables ($z$ in
(\ref{dae2})), whereas the evolution
of the synaptic weights (corresponding to the $y$-variables
in (\ref{dae2})) characterizes the learning scheme; cf.\ \cite{baldi, nnspp}. %Tang y Siljak 
Beyond (\ref{dae2}), constraints may arise not involving the algebraic variables $z$; this
yields a so-called Hessenberg
DAE 
\begin{subequations} \label{hess2}
\begin{eqnarray} 
y' &= & f(y, z) \\
0 & = & g(y).
\end{eqnarray}
\end{subequations}
Such DAEs may also arise naturally in the presence of constraints involving only
node variables (standing for $y$ in (\ref{hess2})) in models which also use certain 
variables supported on the network edges ($z$ in (\ref{hess2})). 
Along the lines later discussed in Section \ref{sec-groups}, Hessenberg DAEs may naturally arise also 
in multilevel networks in which certain restrictions among group variables are not visible 
at the agent's level, and viceversa.
On the other hand, as it happens in nonlinear circuit theory (cf.\ \cite{LMTbook, wsbook}), 
model reduction techniques in general networks
may also lead to quasilinear DAEs of the form $A(x)x'=f(x)$, where $A(x)$
is an everywhere singular matrix-valued map. A reduction to
an explicit ODE formulation may be difficult or impossible in these broader contexts. 
Besides circuit theory, 
where DAEs are now pervasive, different fields are starting to benefit from
the systematic use of constrained formulations
to describe network dynamics; examples are flow networks  \cite{lennartNetworks} 
(including flow dynamics in tree networks %(ok)  
\cite{mayes2011}), %{mayes2011} 
%5) I
integrated process networks %(ok) 
\cite{contou2004}, % {contou2004} 
%7) D
dynamic traffic networks %ok
\cite{han2012},
pressurized water networks \cite{steinbach2005}
%6) Pressurized water networks %(ok)
%\cite{steinbach2005} {steinbach2005} 
or, from a different perspective, 
multiagent descriptor systems \cite{yangliu}.
Analytical results directed to general networks are very promising in these and other related
%research 
areas.

\section{Dynamics on multilevel networks}
\label{sec-groups}

We now turn our attention to a dynamical problem in a
clustered network. As detailed below, we will assume that the set of agents (nodes) is partitioned
into disjoint clusters or groups;
state variables will be defined both on individual agents and on groups, and membership
of a given group will entail a dynamical relation between the node and the group state variables.
Additionally, we will assume 
that there are not only relations among individual agents,
but also among groups.  We are not concerned about the nature of these 
relations or
the criterion according to which the nodes are clustered, but 
if it is of help the reader may think e.g.\ of a set of
economic agents clustered in cooperatives, commercial transactions taking place both among agents
and among cooperatives.
Our goal is to analyze the dynamics of such a two-level network.

\subsection%{Generalized hypergraphs 
{Two-level networks}
\label{subsec-hyper}

%The DNA lady
%Butts
%Clauset (Nature)
%The popes

\noindent We consider a two-level network defined by a 4-tuple $(V, E, H,
E_H),$ where

\begin{itemize}
\renewcommand{\labelitemi}{$\circ$}
%\bullet - default rrr
%\circ — An open circle
%\cdot — A centered dot
%\star — A five-pointed star
%\ast — A centered asterisk
%\rightarrow — A short right-pointing arrow
%\diamondsuit — An open diamond
%For a full list of character commands, see Hypertext Help with LaTeX: Binary and relational operators
\item $V$ is a set of $n \geq 1$ nodes or {\em agents};

\item $E$ is a set of $m$ {\em edges} connecting some pairs of nodes;

\item $H$ is a family of $p \geq 1$ non-empty sets or
{\em groups} of nodes; and

\item $E_H$ is a set of $q$ {\em generalized edges} connecting some 
pairs of groups.
\end{itemize}

\noindent
Note that both $(V, E)$ and $(H, E_H)$ are graphs, whereas $(V, H)$ defines
a hypergraph (with hyperedges corresponding to the above-defined 
groups) if all nodes belong to at least one group
\cite{bergeHyper}. This structure simply models a set of agents 
a) with a dyadic (binary) relation among them; 
b) joined together into certain groups; and
c) displaying also a group-level dyadic relation.
Both $(V, E)$ and $(H, E_H)$ will be assumed to be simple graphs,
and we will focus on cases in which $H$ defines
a {\em partition} of the set of agents %$1, \ldots, n$ 
into $p$ pairwise-disjoint 
groups, the $j$-th one including $n_j \geq 1$ agents. 
Without loss of generality we assume that the agents 
are numbered according to this partition, so that the indices $1, \ldots,
n_1$ correspond to the first group, $n_1+1, \ldots, n_1 + n_2$
to the second one, etc. 

We assume that a reference direction is given to each edge and each generalized edge,
giving both $(V, E)$ and $(H, E_H)$ a digraph structure. 
%*** MAYBE NOT OR NOT HERE: and that both digraphs are connected.  
The incidence matrix describing the dyadic relation at the
agents' level will be denoted by $A$, as in Section \ref{sec-constrained} (cf.\ (\ref{incidence})).
Analogously, the incidence matrix describing the relation among groups will be
written as $A_G=(a_{ij}) \in \R^{p \times q}$, where 
\begin{eqnarray} \label{incidencebis}
a_{ij} = \left\{
\begin{array}{rl}
1 & \text{ if generalized edge } j \text{ leaves  group } i \\
-1 & \text{ if generalized edge } j \text{ enters group } i \\
0 & \text{ if generalized edge } j \text{ is not incident with group } i.
\end{array}
\right.
\end{eqnarray}
In turn, the entries of the incidence matrix $A_H =(a_{ij}) \in \R^{n \times p}$ of the hypergraph $(V, H)$
read as %s as $(a_{ij})$ with
\begin{eqnarray} \label{incidencetris}
a_{ij} = \left\{
\begin{array}{rl}
1 & \text{ if the $i$-th agent belongs to the $j$-th group} \\ 
0 & \text{ otherwise.}
\end{array}
\right.
\end{eqnarray}

\subsection{Node-group dynamics}

In the setting considered above, 
dynamics may take place
both at the agent and at the group level; therefore,
state variables will be defined not only for agents 
but also for groups.
For simplicity, we will assume that each agent and each group
have exactly one (scalar) state variable, to be denoted
by $x_i$ for $i=1, \ldots, n$ and $y_j$ for $j=1, \ldots, p$, respectively.
Alone, the $x$- and the $y$-variables might be understood to
correspond to dynamical processes defined on the graphs 
$(V, E)$ and $(H, E_H)$. Certainly, the interest 
is placed on cases in which there is additionally
a set of dynamic processes relating the
agent state variables $x_i$ with the group ones $y_j$.

In this context, we wish to analyze an extension of
the redistribution scheme discussed in Section \ref{sec-constrained}
to two-level networks. 
Specifically, both at the agent and at the group level
there will be a dynamical process defined by the corresponding Laplacian
matrices, that is, $\ml=A A^{\tra}$ and $\mg=A_G A_G^{\tra}$. 
%(ojo, no $\mh=A_H A_H^{\tra}$).
Without an agent-group dynamic coupling, such dynamical systems would simply read
as 
\begin{eqnarray} \label{consensus-nodes}
x' = -\ml x
\end{eqnarray}
and
\begin{eqnarray} \label{consensus-groups}
y'= -\mg y,
\end{eqnarray}
as in (\ref{consensus}).

The coupling between agents and groups will be defined by the assumption
that  the group variable 
$y_j$ stands for a collectively
stored commodity at the $j$-th group, which 
(in the absence of transactions among agents or groups)
evolves according to the individual amounts $x_i$ stored by the agents
which are clustered in that group. Specifically, we will assume that 
the rate at which the $i$-th agent contributes to $y_j$ (or collects from
$y_j$) is
proportional to the difference 
$$\frac{y_j}{n_j}-x_i.$$
The idea here is that the $i$-th agent contributes to (or collects from) $y_j$ depending
on the difference between its own amount of the resource, $x_i$, and the part  $y_j/n_j$
that would correspond to each agent in an eventual 
uniform distribution of
the group stored commodity $y_j$ to all $n_j$ agents.
% so that asymptotically all agents in the group would tend to have
%the same $x_i$ and the collectively stored commodity would equal 
%the sum of $x$'s.
This yields
$$x_i'=-x_i + \frac{y_j}{n_j}, \text{ for } i \text{ in group } j$$ 
and
$$y_j' = \sum_{i \text{ in } j} x_i - y_j, \ j = 1, \ldots, p.$$
The latter relation is derived immediately from the identity 
$$y_j'=- \sum_{i \text{ in } j} x_i',$$
where for the moment we are disregarding the flows among agents and among groups.

Note that, by construction, 
%For later use notice that the product
$A_H^{\tra}A_H$ is a diagonal matrix of order $p$, the $j$-th diagonal
entry being $n_j$, that is, the number of agents in the $j$-th group.
The
node-group dynamic relations described above can be then globally expressed as
\begin{subequations} \label{node-group}
\begin{eqnarray} 
x' & = &-x + A_H (A_H^{\tra}A_H)^{-1}y \\
y' & = & A_H^{\tra}x - y,
\end{eqnarray}
\end{subequations}
which can be understood as a dynamical process in the hypergraph $(V, H)$.

%\noindent {\bf Node and group dynamics}

The whole dynamical process combines the flows among agents and among groups, described
by (\ref{consensus-nodes}) and (\ref{consensus-groups}), respectively, with the
node-group dynamic relations defined by (\ref{node-group}). This leads to the model
\begin{subequations} \label{full-node-group}
\begin{eqnarray}
x' & = & -(\mathbb{I}_n+\ml)x + A_H (A_H^{\tra}A_H)^{-1}y \\
y' & = & A_H^{\tra}x - (\mathbb{I}_p+\mg)y,
\end{eqnarray}
\end{subequations}
where $\mathbb{I}_n$ and $\mathbb{I}_p$ are identity matrices of orders $n$ and $p$.

\subsection{A stability problem}

System (\ref{full-node-group}) describes a dynamical process in a two-level network,
without any restriction so far in the topology of the graphs $(V, E)$ and $(H, E_H)$.
We address here a stability problem in a simplified setting; we assume that 
there are no interactions at the agents' level or, in graph-theoretic terms, that all
nodes are disconnected. The flow (\ref{consensus-groups}) describing the interactions
among groups in the graph $(H, E_H)$
and the
node-group dynamics (\ref{node-group}) in the hypergraph $(V, H)$
still apply, to yield the model:
\begin{subequations} \label{groupdynamics}
\begin{eqnarray}
x' & = & -x + A_H (A_H^{\tra}A_H)^{-1}y \\
y' & = & A_H^{\tra}x - (\mathbb{I}+\mg)y.
\end{eqnarray}
\end{subequations}
Below we denote by $\chi(i)$ the group to which agent $i$ belongs. 

\begin{theor}\label{th-groups}
Consider system (\ref{groupdynamics})  and assume 
 that the graph $(H, E_H)$, describing the relations among groups,
is connected.

\vspace{3mm}

\noindent {\rm (a)} System (\ref{groupdynamics}) 
has a line of equilibria defined by 
the relations 
\begin{equation} \label{equili}
x_i = \frac{y_{\chi(i)}}{n_{\chi(i)}} \ (i= 1, \ldots, n), \ \
y_1= y_2 =\ldots = y_p.
\end{equation}
%\in \ker A_G$, 
%$x= A_H (A_H^{\tra}A_H)^{-1}y$.

\noindent {\rm (b)} An initial condition $(x(0), y(0))$ converges exponentially to the equilibrium
defined by
\begin{equation} \label{soluti}
x_i = \frac{y_{\chi(i)}}{n_{\chi(i)}}, \ i= 1, \ldots, n, \ y_j = \frac{1}{2p}\left(\dsp\sum_{i=1}^{n}x_i(0) + \dsp\sum_{k=1}^p y_k(0)\right), \ j=1, \ldots, p.
\end{equation}
\end{theor}

\noindent {\bf Proof.} 

\

\noindent (a) By means of a Schur reduction
of (\ref{groupdynamics}), equilibria are
easily checked to satisfy
$$A_H^{\tra}A_H (A_H^{\tra}A_H)^{-1}y - (\mathbb{I}+\mg)y=0,$$
that is,
\begin{equation} \mg y=0. \nonumber
\end{equation}
Since the graph $(H, E_H)$ is connected, $\mg =A_G A_G^{\tra}$ is 
rank-deficient by one, and $\ke \mg = \ke A_G^{\tra}$ is defined by
the relations 
$y_1= y_2 =\ldots = y_p.$ The $x$-components are defined by
$x= A_H (A_H^{\tra}A_H)^{-1}y$. 
By construction $A_H^{\tra}A_H$ is a diagonal matrix, and the
$j$-th diagonal equals $n_j$, that is, the number of
agents in the $j$-th group; we then have
\begin{equation} \label{AHetc}
A_H = \left( \begin{array}{cccccc}
1 & 0 & 0\ \ \ & 0 &\ldots & 0 \\
\vdots & 0 & 0\ \ \  & 0 &\ldots & 0 \\
\text{{\tiny($n_1$)}} & 0 & 0\ \ \  & 0 &\ldots & 0 \\
\vdots & 0 & 0\ \ \  & 0 &\ldots & 0 \\
1 & 0 & 0\ \ \  & 0 &\ldots & 0 \\
0 & 1 & 0\ \ \   & 0 &\ldots & 0 \\
0 &\vdots & 0\ \ \  & 0 & \ldots & 0 \\
0 & \ \text{{\tiny($n_2$)}}\ \ & 0\ \ \  &0 & \ldots & 0 \\
0 &\vdots & 0\ \ \  & 0 & \ldots & 0 \\
0 & 1 & 0\ \ \  & 0 & \ldots & 0 \\
0 & 0 & 1\ \ \  & 0 & \ldots & 0 \\
\vdots & \vdots &\vdots\ \ \   & \vdots & \ddots & \vdots \\
0 & 0 & 0\ \ \ & 0 &\ldots & 1 
\end{array}\right), \  \
(A_H^{\tra}A_H)^{-1} =
\left( \begin{array}{cccccc}
\frac{1}{n_1}\vspace{2mm} & 0 & 0 &\ldots & 0 \\
0 & \frac{1}{n_2}  & 0 &\ldots & 0\vspace{2mm} \\
0 & 0 & \frac{1}{n_3} &\ldots & 0\vspace{2mm} \\
\vdots & \vdots  &\vdots & \ddots  & \vdots\vspace{2mm} \\
0 & 0  & 0 &\ldots & \frac{1}{n_p} \\
\end{array}\right) \ \
\end{equation}
so that
$x= A_H (A_H^{\tra}A_H)^{-1}y$
amounts to 
$x_i = \dsp\frac{y_{\chi(i)}}{n_{\chi(i)}}$ for $i= 1, \ldots, n$,
as stated in (\ref{equili}).

\

\noindent (b) System (\ref{groupdynamics}) will be
proved stable by using Ger\v{s}gorin
theorem \cite{horn}, according to which all eigenvalues of the coefficient
matrix of the right-hand side of (\ref{groupdynamics}), namely
\begin{equation*}
M = \left( \begin{array}{cc}
 -\mathbb{I}  & \hspace{2mm} A_H (A_H^{\tra}A_H)^{-1} \vspace{2mm} \\
 A_H^{\tra} & \hspace{-1mm} - (\mathbb{I}+\mg)
\end{array}
\right),
\end{equation*}
must lie on the union of the discs $|z-m_{jj}| \leq R_j$, with
$$R_j= \sum_{\substack{i=1 \\ i \neq j}}^{n+p} |m_{ij}|.$$
We are denoting by $m_{ij}$ the entries of the matrix $M$. It is obvious that
$$m_{jj}=-1 \text{ if } 1 \leq j \leq n,$$ and 
$$m_{jj}=-1-d_{j-n} \text{ if } n +1 \leq j \leq n+p,$$ 
where $d_j$
stands for the number of connections of the $j$-th group to other groups (that is,
its degree in the graph $(H, E_H)$); we are
making use of the fact that $\mg$ is the Laplacian matrix of the graph $(H, E_H)$ describing the
group connections. 

On the other hand, %it is easy to check that 
there is only
one non-vanishing entry (with value 1) in any column of $A_H^{\tra}$, so that 
$$R_j=1 \text{ if } 1 \leq j \leq n.$$ 
This means that $n$ discs are centered at $-1$ and have radius $1$.
Additionally,  by construction one can see 
that the sum of the entries 
in any column of $A_H (A_H^{\tra}A_H)^{-1}$ is 1; indeed, from (\ref{AHetc}) we may check that 
the $j$-th column of $A_H (A_H^{\tra}A_H)^{-1}$  has $n_j$ nonvanishing
entries equal to $1/n_j$.
Together with
the fact that the off-diagonal entries 
of $-\mg$ define the adjacency matrix of the graph $(H, E_H)$, this yields
$$R_j=1+d_{j-n} \text{ for } n +1 \leq j \leq n+p$$ so that
the remaining $p$ discs  are centered at $-1-d_{j-n}$ and have radius $1+d_{j-n}$ (with $j=n+1, \ldots, n+p$).
For later use, notice that this also shows
that the sum of the entries of any
column of $M$ does vanish.

The remarks above prove that 
all Ger\v{s}gorin discs are located on the left half complex plane,
except for a tangency with the imaginary axis at the origin. From (a) it follows that 
the matrix $M$ has a unique zero eigenvalue (with eigenvectors in $\ke \mg$), 
so that the other $n+p-1$ ones are
away from the imaginary axis. 

The last remark needed to prove (b) is that the quantity 
$\sum_{i=1}^{n}x_i + \sum_{j=1}^p y_j$ is preserved along trajectories,
which is a consequence of the aforementioned fact that the sum of
the entries in each column of $M$ is null (so that $\sum_{i=1}^{n}x_i' + \sum_{j=1}^p y_j'=0$). This means that 
$$\dsp\sum_{i=1}^{n}x_i + \dsp\sum_{j=1}^p y_j=\dsp\sum_{i=1}^{n}x_i(0) + \dsp\sum_{j=1}^p y_j(0)$$ 
is a family of invariant hyperplanes, the evolution rate in each one of them being
characterized by the $n+p-1$ eigenvalues which are located on the left half-plane.
Therefore, the dynamics evolves exponentially towards the point defining
the intersection of such invariant planes and the line of equilibria arising in (a).
By combining the  identities 
$y_1= y_2 =\ldots = y_p$ with the fact that the conditions 
$x_i = y_{\chi(i)}/n_{\chi(i)}$ yield 
$\sum_{i=1}^{n}x_i =\sum_{j=1}^p y_j$ at equilibrium,
it is easy to check that this intersection is given by (\ref{soluti}) and the proof is complete.

\hfill $\Box$

\

\noindent
Theorem \ref{th-groups} extends the results discussed in Section \ref{sec-constrained} to 
the two-level dynamic network here considered. Now the scheme converges asymptotically to an
equal distribution of the resource among all groups, and to a uniform distribution
among the individual agents in each group. Differences may arise between agents of different
groups. 

We leave for future study the full analysis of (\ref{full-node-group}), which may display
more intricate dynamics; further research
may extend these results to problems with flowrates different from the ones yielding the
graph Laplacian matrices in (\ref{consensus-nodes}) and (\ref{consensus-groups}), 
or with node-group dynamic relations different from (\ref{node-group}). 
Our results should also be of interest in an eventual 
deterministic, continuous-time dynamic analysis of 
multilevel networks in general, possibly accommodating
more aggregation levels.

\section{On the equilibrium set of potential-driven flow networks}
\label{sec-eq}

Theorems \ref{th-constrained} and \ref{th-groups} above share a somewhat unusual property
in dynamical systems theory, namely the existence of a line of equilibria. 
Non-isolated equilibria
have received attention within the theory of bifurcation without parameters  
(cf.\ \cite{liebscherbook} and references therein; see also \cite{aulbach}). In our
context, the key element supporting
the existence of a line of equilibria is the fact that systems (\ref{consDAE}) and
(\ref{groupdynamics}) essentially describe the dynamics of potential-driven flows.
As detailed later, equilibria define a one-dimensional manifold (a line) because
of the implicit assumption that the networks involved
are positively weighted. Without this assumption, the structure of the equilibrium
set may be more intricate; in this section we examine
such structure in flow networks 
without this positiveness assumption, arriving at a graph-theoretic characterization
of the structure of the equilibrium set in terms of the network spanning trees. 
In particular, we will apply our results to 
networks based on signed graphs, originally introduced by Harary \cite{harary53}.

\subsection{Nonlinear flows}
\label{subsec-nonlin}

A potential-driven flow on a (directed) graph with $n$ nodes and $m$ edges
is defined by a dynamical system of the form
\begin{eqnarray} \label{flow0}
\alpha_i x_i' = -A_if(A^{\tra}x), \ i=1, \ldots, n.
\end{eqnarray}
Here $x=(x_1, \ldots, x_n)$ and $A_i$ is the $i$-th row of the incidence matrix $A$ as defined in 
(\ref{incidence}). In turn $f:  \R^m\to\R^m$ is a possibly nonlinear, 
differentiable
map describing the flowrates in the edges. It is assumed to have a diagonal structure, that is
$$f=f_1 \times \ldots \times f_m,$$ with $f_j:\R \to \R$ depending only on $(A^{\tra})_jx$,
where $(A^{\tra})_j$ is the $j$-th row of $A^{\tra}$; this means
that the flowrate in edge $j$ (connecting nodes $i$ and $k$) 
depends only on the difference $x_i-x_k$. The variables $x$ can be thought of as a potential 
and hence the ``potential-driven'' label for the flow: electric
potential or pressure are examples in electrical circuits and
water networks, respectively.
Finally, the coefficients $\alpha_i$ are $0$ or $1$ for each node; the case $\alpha_i=1$ (resp.\
$\alpha_i=0$)
corresponds to nodes which accumulate (resp.\ do not accumulate) a certain amount of the 
quantity or resource (e.g.\  
electrical charge, water, gas, a given commodity, etc.) which is flowing in the network;
note that (\ref{flow0}) describes the continuity equations at both types of nodes.

A simple example with $\alpha_i=1$ for all nodes, $f$ being the identity map, is defined
by (\ref{unfold-lapla})-(\ref{unfold-laplb}), yielding the Laplacian dynamics (\ref{consensus}).
A problem combining nodes with $\alpha_i=1$ and $\alpha_i=0$ is 
given by system (\ref{consDAE}), analyzed in Section \ref{sec-constrained}. 
Additionally, the reader may think of 
a (possibly nonlinear) resistive circuit as an example with $\alpha_i=0$ for all nodes.

\subsection{Equilibria and the subimmersion theorem}
\label{subsec-eq}

To make things simpler, in the sequel we assume that $\alpha_i=1$ for all nodes, focusing on 
dynamical systems
of the form
\begin{eqnarray} \label{flow}
x' = -Af(A^{\tra}x), 
\end{eqnarray}
the aforementioned restrictions on the form of $f$ still holding. Also for the sake
of simplicity we assume that the digraph is connected, so that $\rk A=n-1$. For notational
convenience we denote by $F(x)$ the right-hand side of (\ref{flow}), that is,
\begin{eqnarray} \label{F}
F(x) = -Af(A^{\tra}x).
\end{eqnarray}

It is easy to check
that equilibria of (\ref{flow}) may never be isolated; indeed, provided that
%the right-hand side of (\ref{flow}) 
$F(x)$ vanishes for a given $x^*$, 
then $x=x^* + v$ also annihilates $F(x)=-Af(A^{\tra}x)$, 
for any $v \in \ke A^{\tra}$. Note that this kernel is never trivial, being one-dimensional in a
connected digraph.
%its dimension being given by the number of connected components of the digraph.
%Focusing for simplicity on connected digraphs, 
The problem we address in this section is the
characterization of the cases in which the equilibrium set is locally a line,
as it happens in Theorems \ref{th-constrained} and \ref{th-groups}. We show below
that the one-dimensional nature of the equilibrium
set in these theorems is implicitly supported on a positiveness assumption on the 
digraph weights, without which the equilibrium set may locally have a higher dimension. In general,
we do not require $f$ to be linear;
when $f$ is a linear map the equilibrium set is a linear manifold and the results hold
globally.

We will make use of the subimmersion theorem, which states that
if $\Omega$ is an open subset of $\R^n$, and $F$ is a smooth mapping
$\Omega \to \R^p$ such that the Jacobian
matrix $F'(x)$ has constant rank $r \leq p$ on $\Omega$,
then for every $y$ in $F(\Omega)$ the set $F^{-1}(y)$ 
is a submanifold of $\Omega$ with dimension $n-r$ 
(see e.g.\ Th.\ 3.5.17 in \cite{abraham} or Th.\ III.5.8 in \cite{boothby}). %/codimension $d$
%A proof can be found e.g.\ in 
The result also holds if $F$ has constant rank $r$
on a neighborhood of $F^{-1}(y)$.
We will use this result with $F$
given in (\ref{F}), $y=0$, $p=n$ and $r=n-1$
to characterize the problems in which  the equilibrium set
is locally a line around a given equilibrium point $x^*$. The Jacobian matrix $F'(x)$ reads
as
\begin{eqnarray} \label{Fder}
F'(x) = -Af'(A^{\tra}x)A^{\tra}
\end{eqnarray}
and, since $\rk A=n-1$, it follows that $F'(x)$ is persistently rank deficient.
Therefore, for the equilibrium set of (\ref{flow}) 
to be  one-dimensional (at least locally around $x^*$), it is
enough to derive conditions guaranteeing $\rk F'(x^*)=n-1$, since
this maximum rank would  necessarily be attained also on a neighborhood of $x^*$.
In order to examine the rank of $F'(x^*)$, let us denote the derivatives of the
components of $f$ at $x^*$ as
\begin{eqnarray} \label{weights}W_j = f_j'(A^{\tra}x^*), \ j=1, \ldots, m,\end{eqnarray}
and let $W$ stand for the diagonal matrix with entry $W_j$ in the 
$j$-th diagonal position. Note that this matrix is indeed diagonal
because of the assumption that $f_j$ depends only on $(A^{\tra})_jx$.
With this notation we have
\begin{eqnarray} \label{Fder1}
F'(x^*) = -A\W A^{\tra}.
\end{eqnarray}
This expression shows that the Jacobian matrix $F'(x^*)$ has the structure of a 
weighted Laplacian matrix, with weights being defined by the derivatives $f_j'(A^{\tra}x^*)$.

\subsection{Positive weights}
\label{subsec-positive}

If all weights $W_j$ (that is, all derivatives $f_j'(A^{\tra}x^*)$) are positive, then
%$F'(x^*)$ in (\ref{Fder1}) is a negative semidefinite matrix, and 
it is a simple matter to check that
\begin{equation} \ke A\W A^{\tra} = \ke A^{\tra}. \label{kercond}
\end{equation}
Indeed, just note that 
$ A\W A^{\tra}u=0$ implies $u^{\tra}A\W A^{\tra}u=0$ and therefore $A^{\tra}u=0$ because of the
%semidefinite nature of $A\W A^{\tra}$
positiveness of the diagonal matrix $W$. 
The relation depicted in (\ref{kercond}) implies
that $$\rk F'(x^*) = \rk A\W A^{\tra} = \rk A^{\tra} = n-1$$ and, because of the subimmersion
theorem, it follows that equilibria actually define a curve near $x^*$.

Implicitly, this underlies the existence of a line
of equilibria in the setting of Theorems \ref{th-constrained} and \ref{th-groups}.
Note that system (\ref{consensus}), defined by the Laplacian matrix $AA^{\tra}$, can
be understood as a particular instance of the product $A\W A^{\tra}$ with $W=\mathbb{I}_m$.
Actually, one may show that the results of Sections \ref{sec-constrained} and \ref{sec-groups}
actually hold if the product $AA^{\tra}$ is replaced by $A\W A^{\tra}$ with a positive, diagonal
$W$; in other words, those results still apply if the networks are assumed to be positively weighted.
Note that the linear setting considered in those sections yields a 
linear manifold of equilibria and avoids the need for a local approach.

\subsection{Negative weights and the structure of the equilibrium set}

If $W$ in (\ref{Fder1}) includes negative entries, the remarks just stated do not apply, and we
may actually find potential-driven flow dynamics on connected graphs
displaying higher ($\geq 2$) dimensional manifolds of equilibria. 
This may be the case even in linear problems:
examples in a linear context, involving signed graphs, can be found in subsection
\ref{subsec-signed} below. 

A different approach is needed to analyze the structure
of the equilibrium set in digraphs with possibly negative weights.
In this setting, we will drive the problem to a context known in the framework of circuit
theory, namely the analysis of 
nodal admittance matrices, which can be performed along the topological approach stemming from the 
work of J. C. Maxwell \cite{chen97, recski}. 
%Because of  Lemma \ref{lema-red}
In Theorem \ref{th-equilibrium}, which is the main result of this section,
we make use of the notion of the weight of a spanning tree; this is simply
the product of the weights of the tree edges.

%A first step in the characterization of problems with curves of equilibria is the following.

\begin{theor} \label{th-equilibrium}
Assume that the dynamical system (\ref{flow}) is defined on a connected digraph.
The Jacobian matrix $F'(x^*)$ in (\ref{Fder1}) has corank one if and only if
the sum of the spanning tree weights does not vanish, where the edge weights are
given by (\ref{weights}).
If this sum is not null, then the set of equilibria is a curve locally around the equilibrium point $x^*$.
\end{theor}

\noindent {\bf Proof.} %\begin{lema} \label{lema-red}
The matrix  $A \W A^{\tra}$ in the right-hand side of
(\ref{Fder1}) is a weighted Laplacian matrix, which is known to be
rank deficient since $\rk A=n-1$. 
The first step in the proof shows that 
%Let $A \W A^{\tra} \in \R^{n \times n}$ be a weighted Laplacian matrix of a connected digraph, with
%incidence matrix $A \in \R^{n \times m}$ and a diagonal weight matrix $W \in \R^{m \times m}$.
%Then $\rk A \W A^{\tra} < n$, and
all $(n-1)$-minors of $A \W A^{\tra}$ do vanish if and only if a single one
of them does; therefore, the identity $\rk A \W A^{\tra} = n-1$ 
will hold if and only if one (hence any) principal minor
does not vanish.
%\end{lema}

%\noindent {\bf Proof.} The property $\rk A \W A^{\tra} < n$ follows immediately from the 
%fact that $\rk A < n$ for any digraph (actually $\rk A = n-1$ for a connected digraph).
By construction, minors %Horn: determinant of square submatrices
of $A \W A^{\tra}$  of order $n-1$ are determinants of products of the form
$$A_{\text{r}_1} W A_{\text{r}_2}^{\tra},$$
where $A_{\text{r}_1}$ and $A_{\text{r}_2}$ are {\em reduced} incidence matrices, defined by two
arbitrary choices of $(n-1)$ rows of $A$. The key remark 
here comes from a graph-theoretic property, saying that any set of
$n-1$ rows of the incidence matrix of a connected digraph are linearly independent
(actually defining a basis of the so-called {\em cut space}; cf.\ \cite{bollobas}). This means
that a relation of the form 
$$A_{\text{r}_2} = K A_{\text{r}_1}$$ 
holds for a non-singular matrix $K$. Therefore
$$\det (A_{\text{r}_1} W A_{\text{r}_2}^{\tra}) = \det (A_{\text{r}_1} W A_{\text{r}_1}^{\tra}) \hspace{1mm} \det K$$
and it follows that $\det (A_{\text{r}_1} W A_{\text{r}_2}^{\tra})$
does vanish if and only if $\det (A_{\text{r}_1} W A_{\text{r}_1}^{\tra})$ does. Hence,
the eventual vanishing of all $(n-1)$-minors occurs if and only if a single one of them is null;
equivalently, in order to examine
the condition $\rk A \W A^{\tra} = n-1$  it suffices to study the vanishing of a principal minor,
that is, 
the condition %of the form 
%the determinant of an
%arbitrary submatrix of the form \begin{eqnarray}\label{nodal}
%A_{\text{r}_1} W A_{\text{r}_1}^{\tra}, \end{eqnarray}
%as claimed.
%%Because of Lemma \ref{lema-red}, 
%we know that $\rk F'(x^*) = \rk A \W A^{\tra} = n-1$ if and only if
%Our goal is then to examine the condition
\begin{eqnarray}\label{nodalbis}
\det(A_{\text{r}} \hspace{0.5mm}W A_{\text{r}}^{\tra}) \neq 0, \end{eqnarray}
for an arbitrary choice of a reduced incidence matrix $A_{\text{r}}$. Without loss of generality, in what
follows we are allowed to work e.g.\ with the reduced incidence matrix
$A_{\text{r}}$ defined by the first $n-1$ rows of $A$. 
The corank-one condition on $F'(x^*)$
then amounts to evaluating the condition
%non-vanishing of the determinant arising in 
(\ref{nodalbis}).

This can be performed using a Cauchy-Binet determinantal expansion of (\ref{nodalbis}) (cf.\ \cite{horn}), which yields
\begin{equation} \label{expansion}
\det(A_{\text{r}} \hspace{0.5mm}W A_{\text{r}}^{\tra}) = \sum_{\alpha,\beta}\det (A_{\text{r}}^{\alpha}) \det (W_{\alpha}^{\beta}) \det ((A_{\text{r}}^{\tra})_{\beta}),
\end{equation}
the sum being taken over all %ordered 
index sets $\alpha$,
$\beta \subseteq \{1,\dots,m\}$ with $n-1$ elements.  In the submatrices involved
in this expansion, the subscripts (resp.\ superscripts) $\alpha$, $\beta$ 
are used to select a set of rows (resp.\ columns);
for instance, $A_{\text{r}}^{\alpha}$ denotes the submatrix of $A_{\text{r}}$ specified by all rows and
the columns indexed
by $\alpha$, whereas $W_{\alpha}^{\beta}$ is the submatrix of $W$ defined by the rows and columns specified
by $\alpha$ and $\beta$, respectively.

The expansion depicted in (\ref{expansion}) can be simplified using the following two remarks. First,
a set of columns of the reduced incidence matrix $A_{\text{r}}$
is known to yield a non-singular matrix if and only
if such columns correspond to the edges of a spanning tree; in that case, we have 
$\det (A_{\text{r}}^{\alpha})=\pm 1$, because incidence
matrices are totally unimodular  (see e.g.\ \cite{andras1}). 

Second, because of the diagonal form of $W$,
one can easily check that $\det (W_{\alpha}^{\beta})$ does vanish if $\alpha \neq \beta$; when
$\alpha=\beta$, we have
$$\det (W_{\alpha}^{\alpha}) = \prod_{j \in \alpha}W_j,$$
which is the tree weight when $\alpha$ specifies a spanning tree.
Note additionally that, if $\alpha =  \beta$ then $\det (A_{\text{r}}^{\alpha}) = \pm 1$ equals
$\det ((A_{\text{r}}^{\tra})_{\beta})=\det ((A_{\text{r}}^{\tra})_{\alpha})$, so that 
$$\det (A_{\text{r}}^{\alpha}) \det ((A_{\text{r}}^{\tra})_{\alpha})=1.$$

Altogether, these remarks prove that the determinantal
expansion displayed in (\ref{expansion}) amounts to 
$$\det(A_{\text{r}} \hspace{0.5mm}W A_{\text{r}}^{\tra}) = \sum_{\alpha} \prod_{j \in \alpha}W_j,$$
where the sum is taken over the index sets $\alpha$ which specify a spanning tree. 
It follows that %In light of ****************Lemma \ref{lema-red} and 
%the condition depicted in (\ref{nodalbis}), this proves that 
any $(n-1)$-minor of
the Jacobian matrix $F'(x^*)$ is non-null,
%is non-singular, 
and therefore $F'(x^*)$ has corank one, if and only if
the sum of the spanning tree weights does not vanish, as claimed.

Provided that the non-vanishing condition (\ref{nodalbis}) holds, the %resulting 
identity 
$\rk F'(x^*)=n-1$ yields, as a direct consequence of the
subimmersion theorem in the terms stated above, a
local one-dimensional structure for the equilibrium set near $x^*$, and the proof is complete.

\hfill $\Box$

\

\noindent Theorem \ref{th-equilibrium} provides a graph-theoretic characterization, in
terms of the digraph tree structure, 
of the problems which systematically lead to curves (and not higher dimensional or singular
manifolds) of equilibria in potential-driven flow networks with (possibly) some negative weights.
Stability aspects in this context are in the scope of future research; note in particular
that the failing of the non-vanishing requirement in the sum of tree weights might
be responsible for bifurcation phenomena in nonlinear problems.

Certainly, if all weights are positive then
the sum arising in Theorem  \ref{th-equilibrium} is positive and therefore non-null,
because all tree weights are positive. This is of course consistent with
the discussion of subsection \ref{subsec-positive} regarding the one-dimensional nature of
the equilibrium set in positively weighted networks, for which there is no need to use
these tree-based tools.

\subsection{Signed graphs}
\label{subsec-signed}

A nice corollary of Theorem \ref{th-equilibrium} holds for signed graphs, originally introduced
by Harary \cite{harary53} and widely used in balance theory 
%and in the positional analysis of networks 
(cf.\ \cite{harary65} and more recent references such as
\cite{brandes, easley}). A {\em signed graph} or {\em s-graph} is
a graph $(V, E)$ endowed with a map $E \to \{-1, 1\}$, 
that is, an assignment of  
either a $+1$ or a $-1$ weight to all $m$ edges. As above, we will let $W \in \R^{m \times m}$ stand for 
the diagonal matrix of weights.

In this context, a flow problem which arises as a particular case of (\ref{flow}) is
\begin{eqnarray} \label{flowsigned}
x' = -A\W A^{\tra}x.
\end{eqnarray}
This system can be again understood as a flow, via the continuity equations
$x'=-Au$, in which the flowrates $u$ are given
by the linear relation $\W A^{\tra}x$. With respect to (\ref{unfold-lapla})-(\ref{unfold-laplb}), yielding
the Laplacian dynamics (\ref{consensus}), now the presence of $-1$ weight values models
a flow in which adjacent agents $i$, $k$ tend to {\em increase} the difference between
$x_i$ and $x_k$, since
the flow from $i$ to $k$ now equals $x_k - x_i$ (instead of $x_i-x_k$, as in subsection \ref{subsec-flows}).

Contrary to the results in Sections \ref{sec-constrained} and \ref{sec-groups}, in which equilibria
define a line, the presence
of negative weights might now result in higher dimensional equilibrium manifolds. The cases in which
this may happen are exactly characterized in Corollary \ref{coro-signed} below. Since all weights
are $+1$ or $-1$, we may now define a spanning tree as {\em positive} or {\em negative} simply if
its weight product is $+1$ or $-1$, respectively, or, equivalently, if it contains an even (resp.\ odd) 
number of edges with negative weight.

\begin{coro} Let the dynamical system (\ref{flowsigned}) be defined on a connected signed
digraph.
%and assume that the entries of the diagonal matrix $W$ are either $+1$ or $-1$. 
Then,
the dimension of the equilibrium set is higher than one if and only if the 
numbers of positive and negative spanning trees coincide.
%number of positive trees %exactly matches 
%equals the number of negative trees.
\label{coro-signed}
\end{coro}

%There is no need for a proof since t
This results follows immediately from Theorem \ref{th-equilibrium}
since, for the sum of 
weight products to vanish in a signed graph, the amount of positive trees (which are responsible for a
$+1$ term in the sum) must obviously match the number of negative trees (which yield a $-1$ term in the sum).

\

\begin{figure}[h]
\vspace{4mm}
\centering
\epsfig{figure=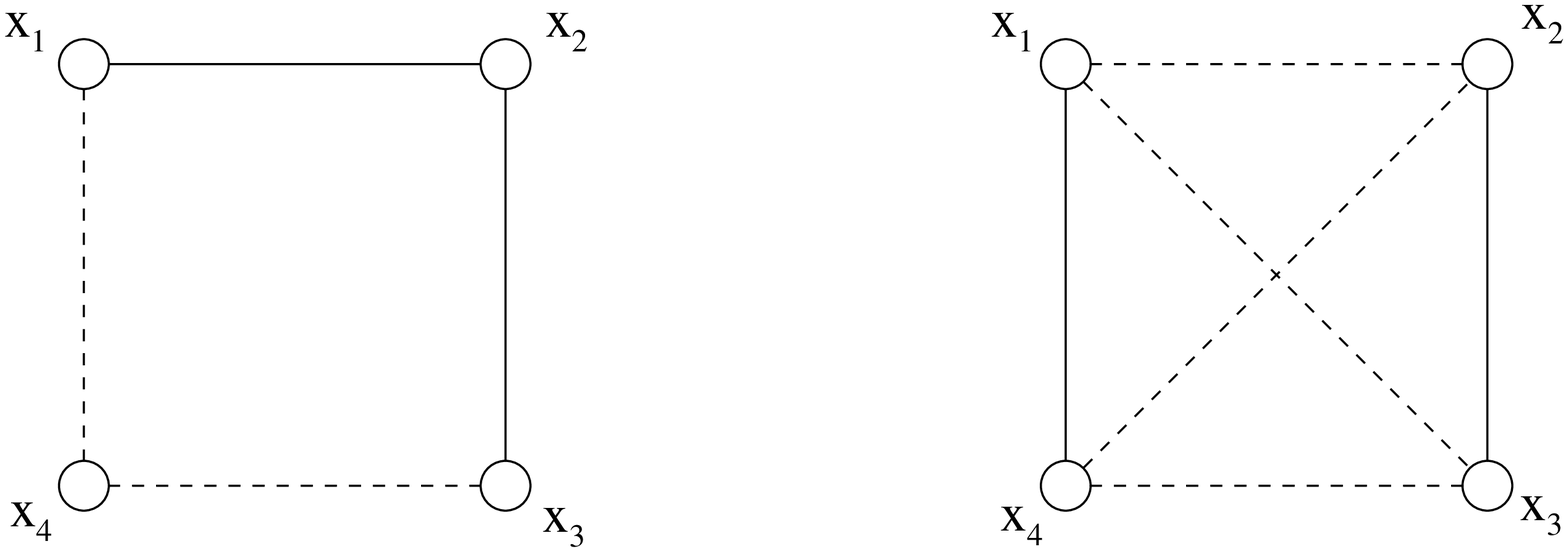, width=0.64\textwidth}%, angle=270}
\vspace{1mm}
\caption{Sign assignments yielding degenerate flow dynamics on (a) a 4-cycle; (b) $K_4$. Edges with a
negative sign are dashed.}% (pdte***)} %defined by the P-structure.}
\label{fig-signed}
\vspace{3mm}
\end{figure}

\noindent{\bf Examples.} Simple examples illustrating the result above can be defined using
the graphs shown in Figure \ref{fig-signed}. The first case, on the left of the figure, 
is simply a 4-cycle with two
positive signs (continuous lines) and two negative signs (dashed lines).
This graph has just four spanning trees, two of which are positive and the other two negative. 
The second example, on the right, depicts a complete graph $K_4$ with two
positive and four negative signs; this is a so-called balanced structure
\cite{harary53}, in which the nodes may be split in two groups in a way such that
all edges inside a group have positive weights and all connections between
groups are negative. In this example the groups are defined by
the two nodes on the left and those on the right, respectively.
According to Cayley's
formula this graph has $4^2=16$ spanning trees, and it is not difficult to check
that exactly half of them are positive.
This
means that the dynamical system (\ref{flowsigned}) should exhibit in both cases
a linear manifold of equilibria with dimension greater than one.

In the first case, we number and orientate each edge beginning on the top and 
according to a clockwise orientation of the cycle. This yields
\begin{equation*} %\label{signed1}
A = \left( \begin{array}{rrrr}
1 & 0 & 0 & -1 \\
-1 & 1 & 0 & 0 \\
0 & -1 & 1 & 0 \\
0 & 0 & -1 & 1
\end{array}
\right).
\end{equation*}
With $W=\mathrm{diag}\hspace{1mm}(1,1,-1,-1)$, 
the right-hand side of  (\ref{flowsigned})
reads in this case as
\begin{equation*} %\label{signed1}
-A \W A^{\tra}=\left( \begin{array}{rrrr}
0 & 1 & 0 & -1 \\
1 & -2 & 1 & 0 \\
0 & 1 & 0 & -1 \\
-1 & 0 & -1 & 2
\end{array}
\right),
\end{equation*}
The kernel of this matrix defines the equilibrium set and is defined by the relations
\begin{equation*} %\label{signed1}
x_2 = x_4 = \frac{x_1 + x_3}{2},
\end{equation*}
%This is a simple example which allows for a straightforward explanation.
hence defining a two-dimensional linear manifold (a plane) of equilibria. Notice that
these equilibrium points arise from a
constant flow $x_i-x_{i+1}$ ($i=1, \ldots, 4$, with the terminological abuse $x_5=x_1$) 
which annihilates all derivatives $x_i'$, 
flowing in the clockwise (resp.\ counterclockwise) direction if $x_1 > x_3$
(resp.\ if $x_1 < x_3$). Note that the equilibrium plane includes the line 
$x_1=x_2=x_3=x_4$ for which the flow vanishes.

In the example defined on $K_4$ (on the right of Figure \ref{fig-signed}), 
we number and orientate the edges according to the
following incidence matrix:
\begin{equation*} %\label{signed1}
A = \left( \begin{array}{rrrrrr}
1 & 0 & 0 & -1 & 1 & 0 \\
-1 & 1 & 0 & 0 & 0 & 1 \\
0 & -1 & 1 & 0 & -1 & 0\\
0 & 0 & -1 & 1 & 0 & -1
\end{array}
\right),
\end{equation*}
and the weights $W=\mathrm{diag}\hspace{1mm}(-1,1,-1,1,-1,-1)$
yield
\begin{equation*} %\label{signed1}
-A \W A^{\tra}=\left( \begin{array}{rrrr}
1 & -1 & -1 & 1 \\
-1 & 1 & 1 & -1 \\
-1 & 1 & 1 & -1 \\
1 & -1 & -1 & 1 
\end{array}
\right).
\end{equation*}
Now the equilibrium set is three-dimensional, being defined by the identity
\begin{equation} 
x_1 + x_4 = x_2 + x_3. \label{K4}
\end{equation}
Again, these equilibrium solutions yield nonvanishing flows in the graph
edges, except for those in the line $x_1=x_2=x_3=x_4$.
This example also illustrates that the rank drop in the
product $A \W A^{\tra}$ may be greater than two. 
%Solutions with nontrivial flows now keep

\

We leave for future work the analysis of the nature of such equilibrium solutions 
exhibiting non-vanishing flows, that is, equilibria which are not in 
the line $x_1= x_2 = \ldots = x_n$; 
this analysis should explain, for instance, the relation 
of the equilibrium solutions in the second example with the presence of a balanced
structure in the graph: 
note that (\ref{K4}) expresses that the total amount of the resource
$x$ stored in both groups (namely, $\{x_1, x_4\}$ and $\{x_2, x_3\}$) 
is the same. More generally, further research might address 
stability aspects and bifurcations in these negatively weighted networks, not
only in signed graphs but specially in the nonlinear setting introduced in subsection
\ref{subsec-nonlin}.

Note finally that even though in Section \ref{sec-eq} we have focused for simplicity on flow networks with
homogeneous agents (because of the assumption $\alpha_i=1$ for all nodes,
cf.\ subsection  \ref{subsec-eq}), the results are also of potential interest in broader contexts,
involving e.g.\ heterogeneous nodes or group dynamics as in Sections \ref{sec-constrained}
and \ref{sec-groups}.

\section{Concluding remarks}
\label{sec-con}

We have addressed in this paper some qualitative problems involving network
dynamics beyond the scenario defined by 
a set of dynamical systems supported on the
nodes of a graph. First, the attention has been focused on potential-driven flow dynamics
either on networks with heterogeneous agents or on multilevel networks. Differential-algebraic
models, arising here in flow
networks with heterogeneous agents, are worth receiving further
attention in order to model and 
analyze constrained dynamics on general networks. Dynamical systems on multilevel networks
also have a broad scope for future research. From a different perspective, we have analyzed
the structure of the equilibrium set in (possibly nonlinear) flow networks with negative weights.
The results in this regard apply in particular to signed graphs, and the connection
of our results to the theory of balance and clusterability in signed graphs requires 
further analysis. In greater generality, stability properties 
and bifurcations are worth being studied in networks with negative weights.

Our research is of potential interest in the analysis of social and economic networks involving
flows. Future work might also extend these results to problems displaying network evolution,
and also to stochastic and/or discrete-time contexts.


\begin{thebibliography}{99}
\bibitem{abraham} R. Abraham, J. E. Marsden and T. Ratiu, 
{\em Manifolds, Tensor Analysis, and Applications},
Springer-Verlag, 1988.

\bibitem{albertbarabasireview} R. Albert and
A.-L. Barab\'asi,
Statistical mechanics of complex networks,
{\em Rev. Modern Physics} {\bf 74} (2002) 47-97.

\bibitem{andras1} B. Andr\'{a}sfai, {\em Introductory Graph Theory,} 
Akad\'{e}miai Kiad\'{o}, Budapest, 1977.

%\bibitem{arkak2007} M. Arcak, Passivity as a design tool for group coordination,
%{\em IEEE Trans. Automatic Control} {\bf 52} (2007) 1380-1390.

\bibitem{aulbach} B. Aulbach, {\em Continuous and Discrete Dynamics near
Manifolds of Equilibria}, Lect. Note Math. 1058, Springer-Verlag, 1984.


\bibitem{baldi} P. Baldi, Gradient descent learning algorithm overview: A
general dynamical systems perspective, {\em IEEE Trans. Neural Networks}
{\bf  6} (1995) 182-195.


\bibitem{barabasi99} A.-L. Barab\'asi and R. Albert,
Emergence of scaling in random networks, {\em Science} {\bf 286}
(1999) 509-512.

\bibitem{barrat} A. Barrat, M. Barthel\'emy and A. Vespignani,
{\em Dynamical Processes on Complex Networks,} 
Cambridge Univ.\ Press, 2008.

\bibitem{bergeHyper} C. Berge, {\em Hypergraphs}, North-Holland, 1989.

\bibitem{bollobas} B. Bollob\'{a}s, {\em Modern Graph Theory},
Springer-Verlag, 1998.


\bibitem{boothby} W. M. Boothby, {\em An Introduction to 
Differentiable Manifolds and Riemannian Geometry},
Academic Press, 1986.

\bibitem{brandes} U. Brandes and T. Erlebach (eds.),
{\em Network Analysis. Methodological Foundations,} Springer-Verlag, 2005.


\bibitem{bre1} K. E. Brenan, S. L. Campbell and L. R. Petzold,
{\em Numerical Solution of Initial-Value Problems in
Differential-Algebraic Equations}, %Classics in Applied Mathematics, 
SIAM, 1996.

\bibitem{mmsna} P. J. Carrington, J. Scott and
S. Wasserman (eds.), {\em Models and Methods in Social Network Analysis},
Cambridge Univ. Press, 2005.

\bibitem{chen97} W. K. Chen, {\em Graph Theory and its Engineering
Applications,} World Scientific, 1997.

%\bibitem{chuamemristor71} L. O. Chua, Memristor -- The missing circuit
%  element, {\em IEEE Trans. Circuit Theory} {\bf 18} (1971) 507-519.

\bibitem{contou2004} M. N. Contou-Carrere and P. Daoutidis,
Dynamic precompensation and output feedback control of
integrated process networks, {\em Proc.\ American Control Conf. 2004}, pp.\
2909-2914, 2004.

\bibitem{easley} D. Easley and J. Kleinberg, {\em Networks, Crows and Markets}, 
Cambridge Univ. Press, 2010.

\bibitem{erdos1959} P. Erd\H{o}s and A. R\'enyi, 
On random graphs I, {\em Publicationes Mathematicae} {\bf 6} (1959) 290-297.

\bibitem{erdos1960} P. Erd\H{o}s and A. R\'enyi, 
On the evolution of random graphs,
{\em Publ. Math. Inst. Hungarian Academy Sci.} {\bf 5} (1960) 17-61.

\bibitem{fronczak} 
P. Fronczak, P., A. Fronczak and J. A. Holyst, Self-
organized criticality and coevolution of network structure
and dynamics, {\em Phys. Rev. E} {\bf 73} (2006), 046117-4. 
%(doi:10.1103/PhysRevE.73.046117)

\bibitem{gross} T. Gross and B. Blasius, 
Adaptive coevolutionary networks: a review,
%http://rsif.royalsocietypublishing.org/content/5/20/259.full
{\em J. R. Soc. Interface} {\bf 5} (2008)
259–271.


\bibitem{han2012} K. Han, B. Piccoli, T. L. Friesz and T. Yao, 
A continuous-time link-based kinematic wave model for dynamic traffic networks,
Preprint 1208.5141, ArXiv, 2012. %http://arxiv.org/abs/1208.5141

\bibitem{harary53} F. Harary, 
On the notion of balance of a signed graph, {\em Michigan Math. J.} {\bf 2} (1953) 143-146.

\bibitem{harary65} F. Harary, R. Z. Norman and D. Cartwright,
{\em Structural Models. An Introduction to the Theory of Directed Graphs,}
John Wiley \& Sons, 1965. 

\bibitem{haykin} S. O. Haykin,
{\em Neural Networks and Learning Machines}, Pearson, 2009.

\bibitem{holmenewman}
P. Holme and M. E. J. Newman, Nonequilibrium phase
transition in the coevolution of networks and opinions,
{\em Phys. Rev. E} {\bf 74} (2007) 056108-5. 
%(doi:10.1103/PhysRevE.74.056108)

\bibitem{horn} R. A. Horn and Ch. R. Johnson,
{\em Topics in Matrix Analysis,} Cambridge University Press, 1991.

\bibitem{lennartNetworks} L. Jansen and C. Tischendorf, 
A unified (P)DAE modeling approach for flow networks, in
S. Sch\"ops, A. Bartel, M. G\"unther, E. J. W. ter Maten and P. C. M\"uller (eds.), 
{\em Progress in Differential-Algebraic Equations,} pp. 127-151, Springer, 2014.

\bibitem{kaob} Y. Kaob and C. Wang, 
Global stability analysis for stochastic coupled reaction–diffusion 
systems on networks, {\em Nonlinear Analysis: Real World Appl.} {\bf 14} (2013) 1457-1465.

\bibitem{koenigbattiston} M. D. K\"onig and S. Battiston,
From graph theory to models of economic networks. A tutorial,
in A. K. Naimzada {\em et al.} (eds.), {\em Networks, Topology and Dynamics},
pp. 23-63, Springer, 2009.

\bibitem{kmbook} P. Kunkel and V. Mehrmann, {\em
  Differential-Algebraic 
Equations. Analysis and Numerical Solution},
EMS, 2006.

\bibitem{LMTbook} R. Lamour, R. M\"arz and C. Tischendorf,
{\em Differential-Algebraic Equations: A Projector Based Analysis,} 
DAE Forum, Springer, 2013.


\bibitem{lieberman} E. Lieberman, C. Hauert and M. A. Nowak,
Evolutionary dynamics on graphs, {\em Nature} {\bf 433} (2005) 312-316.

\bibitem{liebscherbook} S. Liebscher, {\em Bifurcation without Parameters,} Springer, 2015.



\bibitem{liu2011} Y. Y. Liu, J. J. Slotine and A. L. Barab\'asi, 
Controllability of complex networks, {\em Nature}
{\bf 473} 167-173, 2011.

\bibitem{liu2013} Y. Y. Liu, J. J. Slotine and A. L. Barab\'asi, 
Observability of complex systems, {\em Proc.\ Nat.\ Acad. Sciences USA}
{\bf 119} 2460-2465, 2013.

\bibitem{lozano} S. Lozano,
Dynamics of social complex networks: Some insights into recent research,
in {\em  Dynamics On and Of Complex Networks}, Springer, 2009, pp. 133-143.

\bibitem{mayes2011} %draft webmail
J. Mayes and M. Sen, Approximation of potential-driven flow dynamics in
large-scale self-similar tree networks, {\em Proc.\ R.\ Soc.\ A}
{\bf 467} (2011) 2810-2824. 



\bibitem{newmanreview}
 M. E. J. Newman, The structure and function of complex networks,
{\em SIAM Review} {\bf 45} (2003) 167-256.

\bibitem{newmannetworks}
 M. E. J. Newman, {\em Networks: An Introduction,}
Oxford Univ.\ Press, 2010. 

%\bibitem{olfati2003} R. Olfati-Saber and R. M. Murray,
%Consensus protocols for networks of dynamic agents, 
%{\em Proc. Amer. Control Conf. 2003}, pp. 951-956.
%Tienen otro en IEEE_TAC 04

\bibitem{olfati2007} R. Olfati-Saber, J. A. Fax and R. M. Murray,
Consensus and cooperation in networked multi-agent systems,
{\em Proc. IEEE} {\bf 95} (2007) 215-233.

\bibitem{price65}
D. J. S. Price, Networks of scientific papers, {\em Science} {\bf 149}
 (1965) 510-515.

\bibitem{price76} D. J. S. Price, 
A general theory of bibliometric and other cumulative advantage processes,
{\em J. Amer. Soc. Inform. Sci.} {\bf 27} (1976) 292-306.

%\bibitem{rabrheintheo} P. J. Rabier and W. C. Rheinboldt,
%Theoretical and numerical analysis of differen\-tial-algebraic
%equations, {\em Handbook of Numerical Analysis,} Vol. VIII, 
%183-540, North-Holland (2002).

\bibitem{rahmani2009}
A. Rahmani, M. Ji, M. Mesbahi and M. Egerstedt, Controllability of multi-agent systems
from a graph-theoretic perspective, {\em SIAM J. Control Optim.} {\bf 48} (2009) 162-186.

\bibitem{recski} A. Recski, {\em Matroid Theory and its Applications in Electric
 Network
Theory and in Statics}, Springer-Verlag, 1989.



\bibitem{wsbook} R. Riaza, {\em Differential-Algebraic Systems,}
%Analytical Aspects and Circuit Applications}, 
World Scientific, 2008.

\bibitem{nnspp} R. Riaza and P. J. Zufiria, 
Differential-algebraic equations and singular perturbation 
methods in recurrent neural learning, {\em Dynamical Systems} {\bf 18}
(2003) 89-105.

%\bibitem{sikolya} E. Sikolya, Flows in networks with dynamic
%ramification nodes, {\em J. Evol. Equ.} {\bf 5} (2005) 441-463.

\bibitem{siljak} D. D. Siljak, Dynamic graphs,
{\em Nonlinear Analysis: Hybrid Systems} {\bf 2} (2008) 
544-567.


\bibitem{steinbach2005} M. C. Steinbach, Topological index criteria
in DAE for water networks, Preprint 05-49, 
Konrad-Zuse-Zentrum f\"ur Informationstechnik 
Berlin, 2005.

\bibitem{strogatzreview}
S. H. Strogatz, Exploring complex networks, {\em Nature}
{\bf 410}
 (2001) 268–276. %on-line on the web

\bibitem{tanner2004} H. G. Tanner, On the controllability of nearest neighbor interconnections, {\em Proc. IEEE Conf. Decision and Control 2004},
pp. 2467–2472, 2004.

\bibitem{yangliu} X. R. Yang
and G. P. Liu,
Necessary and sufficient consensus conditions of 
descriptor multi-agent systems,
{\em IEEE Trans. Cir. Sys. I} {\bf 59} (2012) 
2669-2677. 

\bibitem{zelazo} D. Zelazo and M. Mesbahi,
Graph-theoretic methods for networked dynamic systems:
Heterogeneity and ${\mathcal H}_2$ performance,
in {\em Efficient Modeling and Control of Large-Scale Systems},
Springer, 2010, pp 219-249. 

\end{thebibliography}
\end{document}